\documentclass[12pt]{article}
\usepackage{theorem,amsfonts,amssymb}
\newtheorem{theorem}{Theorem}
\newtheorem{proposition}{Proposition}
\newtheorem{lemma}{Lemma}
\newtheorem{corollary}{Corollary}
\theorembodyfont{\upshape}
\newtheorem{definition}{Definition}
\newtheorem{remark}{Remark}
\newtheorem{example}{Example}
\newtheorem{proof}{Proof}

\newcommand{\bt}{\begin{theorem}}
\newcommand{\et}{\end{theorem}}
\newcommand{\bl}{\begin{lemma}}
\newcommand{\el}{\end{lemma}}
\newcommand{\bp}{\begin{proposition}}
\newcommand{\ep}{\end{proposition}}
\newcommand{\bex}{\begin{example}}
\newcommand{\eex}{\end{example}}
\newcommand{\bc}{\begin{corollary}}
\newcommand{\ec}{\end{corollary}}
\newcommand{\bo}{\begin{proof}}
\newcommand{\eo}{\end{proof}}
\newcommand{\bd}{\begin{definition}}
\newcommand{\ed}{\end{definition}}
\newcommand{\br}{\begin{remark}}
\newcommand{\er}{\end{remark}}
\newcommand{\be}{\begin{enumerate}}
\newcommand{\ee}{\end{enumerate}}

\begin{document}

\title{Ad~$(G)$ is of type $R$ implies $G$ is of type $R$ for certain $p$-adic Lie groups}
\author{C. R. E. Raja}
\date{}
\maketitle

\let\ol=\overline
\let\epsi=\epsilon
\let\vepsi=\varepsilon
\let\lam=\lambda
\let\Lam=\Lambda
\let\ap=\alpha
\let\vp=\varphi
\let\ra=\rightarrow
\let\Ra=\Rightarrow
\let \Llra=\Longleftrightarrow
\let\Lla=\Longleftarrow
\let\lra=\longrightarrow
\let\Lra=\Longrightarrow
\let\ba=\beta
\let\ov=\overline
\let\ga=\gamma
\let\Ba=\Delta
\let\Ga=\Gamma
\let\da=\delta
\let\Oa=\Omega
\let\Lam=\Lambda
\let\un=\upsilon

\newcommand{\Ker}{{\rm Ker}}
\newcommand{\Ad}{{\rm Ad}}
\newcommand{\cK}{{\cal K}}
\newcommand{\pc}{T_{si}}
\newcommand{\Aut}{{\rm Aut}}
\newcommand{\cR}{{\cal R}}
\newcommand{\Spr}{{\rm Spr}}
\newcommand{\cG}{{\cal G}}
\newcommand{\cH}{{\cal H}}
\newcommand{\G}{{\mathbb G}}
\newcommand{\Z}{{\mathbb Z}}
\newcommand{\Q}{{\mathbb Q}}
\newcommand{\cN}{{\cal N}}
\newcommand{\cS}{{\cal S}}
\newcommand{\N}{{\mathbb N}}
\newcommand{\R}{{\mathbb R}}
\newcommand{\C}{{\mathbb C}}
\newcommand{\T}{{\mathbb T}}
\newcommand{\mP}{{\mathbb P}}

\begin{abstract}
We provide a sufficient condition for a $p$-adic Lie group to be of
type $R$ when its adjoint image is of type $R$.
\end{abstract}

\medskip
\noindent{\it 2000 Mathematics Subject Classification:} 22E15,
22E20.

\medskip
\noindent{\it Key words: $p$-adic Lie group, type $R$, adjoint
representation. }

Let $G$ be a $p$-adic Lie group and $\cG$ be the Lie algebra of $G$.
Then there is an analytic homomorphism ${\rm Ad}\colon G \to GL(\cG
)$, called adjoint representation, such that Ad~$(x)$ is the
differential of the inner-automorphism - $g\mapsto xgx^{-1}$ - on
$G$ defined by $x$.  We refer to \cite{Bo} and \cite{Se} for basics
and results concerning $p$-adic Lie groups

We say that a $p$-adic Lie group $G$ is of type $R$ if the
eigenvalues of Ad~$(x)$ are of $p$-adic absolute value one for any
$x\in G$.

\bex

We now give some examples of $p$-adic Lie groups.

\be

\item [(1)] Abelian groups: $p$-adic vector spaces such as $\Q _p ^k$.

\item [(2)] $p$-adic Heisenberg group: $\{ (a, x, z) \mid
a, x \in \Q _p ^k, ~~z\in \Q _p \}$
with multiplication given by $$(a, x, z)(a', x', z') = (a+a', x+x',
z+z'+<a,x'>)$$ where $<u,v>= \sum u_iv_i$ for any two $u,v\in \Q _p
^k$.

\item [(3)] The solvable group $(\Q _p\setminus \{ 0 \} ) \ltimes \Q
_p$: $\{ (a, x) \mid a \in \Q _p \setminus \{ 0 \}, ~~x\in \Q _p \}$
with multiplication given by $$(a, x)(a', x') = (aa', x+ax').$$

\item [(4)] In general any closed subgroup of $GL_n(\Q_p)$ (see
\cite{Se} and \cite{Bo}).

\ee

\eex

Among the above examples, (1) and (2) are of type $R$ but (3) is not
of type $R$.

For an automorphism $\ap$ of $G$, we define the subgroups $U_{\pm}
(\ap )$ by $U_{\pm} (\ap )= \{ x\in G \mid \lim _{n\to \pm \infty}
\ap ^n (x) = e\}$: when only one automorphism is under
consideration, we write $U_{\pm}$ instead of $U_{\pm }(\ap )$.  The
study of these subgroups for innerautomorphisms plays a crucial role
in proving type $R$. For instance, the following observation which
is essentially contained in \cite{Wa} and Theorem 1 of \cite{Ra1} .

\bp\label{prl} A $p$-adic Lie group $G$ is of type $R$ if and only
if $U_{\pm} (\ap )$ is trivial for any innerautomorphism $\ap$.\ep

\bo If $G$ is of type $R$ and $\ap$ is any innerautomorphism, then
by Theorem 3.5 and Corollary 1 of \cite{Wa} we get that $G$ has
arbitrarily small compact open subgroups stable under $\ap$.  Since
$U_\pm (\ap )$ is contained in any $\ap$-stable open subgroup,
$U_{\pm} (\ap )$ is trivial.

Conversely, suppose $U_\pm (\ap )$ is trivial for any
innerautomorphism $\ap$ of $G$.  Then by Theorem 3.5 (iii) of
\cite{Wa} we see that condition (2) of Theorem 1 in \cite{Ra1} is
satisfied and hence $G$ is of type $R$. \eo

In this note we would like to explore the question: {\it Ad~$(G)$ is
of type $R$ implies $G$ is of type $R$.}

The answer to this question is positive if the kernel of Ad is the
center of $G$ but for a $p$-adic Lie group kernel of Ad need not be
the center of $G$ (see Example \ref{h}).  It may be noted that for
Zariski-connected $p$-adic algebraic groups, kernel of Ad is the
center of $G$.

The following example shows that Ad~$(G)$ is of type $R$ need not
imply $G$ is of type $R$, that is answer to our question is not
always positive.

\bex\label{h}

Let $G= \{ (n, a,x,z+\Z _p) \mid n \in \Z, a, x , z \in \Q _p \}$
with multiplication given by $$(n, a,x,z+\Z _p)(m, b,y,z'+\Z _p) =
(n+m, a+p^nb,x+p^{-n}y,z+z'+p^{-n}<a,y>+ \Z _p).$$ Let $\ap \colon G
\to G$ be the innerautomorphism defined by $(n,0,0,0)$.  Then $U_+ =
\{ (0,a, 0, 0) \mid a \in \Q _p \}$ and $U_-= \{ (0,0, x, 0) \mid x
\in \Q _p \}$. But Ad~$(G) \simeq \Z$, hence Ad~$(G)$ is of type $R$
but $G$ is not of type $R$ as $p$ and $p^{-1}$ are eigenvalues of
innerautomorphism given by $(n,0,0,0)$.

\eex

Following provides a sufficient condition for affirmative answer.

\bt\label{t1} Let $G$ be a $p$-adic Lie group such that Ad~$(G)$ is
of type $R$ and open subgroups of $G$ and their quotients are
unimodular. Suppose any topologically finitely generated subgroup of
the kernel of the adjoint representation is a discrete extension of
a $k$-step solvable group for some $k>0$.  Then $G$ is also of type
$R$. \et

\br The condition that open subgroups of $G$ and their quotients are
unimodular is necessary.  If $G$ is of type $R$, then any open
subgroup and its quotients are also of type $R$.  Now it follows
from Theorem 1 and Corollary 1 of \cite{Ra1} that open subgroups of
$G$ and their quotients are unimodular. \er

We first introduce a few notations.  For any two closed subgroups
$A$ and $B$ of a Hausdorff topological group $X$.  Let $C_0(A, B) =
\langle A, B\rangle$ denote the closed subgroup generated by $A$ and
$B$ in $X$, and $C_1(A, B)= \ol {[A, B]}$, $C_k(A, B) = \ol
{[C_{k-1}(A, B), C_{k-1}(A,B)]}$ for all $k
>1$.  If $A=B$, let $C_k(A)=C_k(A, A)$ for all $k \geq 1$.
$X$ is called a $k$-step solvable group if $C_k(X)=\{ e \}$.

We next develop a few general results.

\bl\label{l1} Let $G$ be a $p$-adic Lie group such that Ad~$(G)$ is
of type $R$.  For $g\in G$, if $\ap$ is the innerautomorphism of $G$
defined by $g$, then there is a $\ap$-invariant compact subgroup $O$
such that  $O$ is centralized by $U_{\pm} (\ap)$ and $U_+(\ap )OU_-
(\ap )$ is an open subgroup of $G$. \el

\bo Let $\ap$ be the inner-automorphism of $G$ defined by $g\in G$
and $U_{\pm} = \{ x\in G \mid \lim _{n\to \pm \infty} \ap ^n (x) =
e\}$. Since Ad~$(G)$ is of type $R$, $U_\pm \subset {\rm Ker~(Ad)}$.
It follows from \cite{Wa} that $U_\pm$ are unipotent closed
subgroups of $G$. Since Lie algebra of ${\rm Ker~(Ad)}$ is the
center of the Lie algebra of $G$, $U_\pm$ is abelian.  Since Lie
algebra of ${\rm Ker~(Ad)}$ is abelian, ${\rm Ker~(Ad)}$ contains a
compact open abelian subgroup (cf. Corollary 3, Section 4.1, Chapter
III of \cite{Bo}). Thus, there exists compact open subgroups $K_\pm
\subset U_\pm$ such that $K_+K_-$ is an abelian group.  Since $\ap
^{\pm 1}|_{U_\pm}$ is a contraction, we may assume that $\ap ^{-\pm
i}(K_\pm) $ is increasing and $U_\pm = \cup _{i\geq 0} \ap ^{-\pm
i}(K_\pm)$. Since $K_\pm $ are compact subgroups of $U_\pm$ which
are $p$-adic vector spaces, $K_\pm$ have a dense finitely generated
subgroup. Hence $K_+K_-$ has a dense finitely generated subgroup.
Since $K_+K_-$ is contained in the kernel of the adjoint
representation, each element of $K_+K_-$ centralize an open subgroup
of $G$ (cf. Theorem 3, Section 7, Chapter III of \cite{Bo}).  Since
$K_+K_-$ has a dense finitely generated subgroup, $K_+K_-$
centralize an open subgroup of $G$.

Let $M= \{ x\in G \mid \ol{\{ \ap ^n(x) \}} ~~{\rm is ~~ compact}
~~\}$.  Then $M$ is a closed $\ap$-stable subgroup of $G$ and $M$
contains arbitrarily small $\ap$-stable compact open subgroups.
Since $K_+K_-$ centralize an open subgroup of $G$, there exists a
$\ap$-stable compact open subgroup $O$ of $M$ such that $O$ is
centralized by $K_+K_-$.  Since $\ap (O)=O$ and $\cup \ap ^n (K_\pm)
= U_\pm$, we get that $O$ is centralized by $U_\pm$.  It also
follows from Theorem 3.5 (iii) of \cite{Wa} that $K_+OK_-$ is an
open (subgroup) in $G$.

\eo

\bl\label{topgp} Let $G$ be a $p$-adic Lie group and $\ap$ be an
automorphism of $G$ such that the closed subgroup $P$ generated by
$U_{\pm}(\ap )$ is a discrete extension of a solvable group. Suppose
Haar measures on $P$ and any of its quotients is $\ap$-invariant.
Then $P$ is trivial. \el

\bo The group $Q:=P/C_1(P)$ is abelian and $\ap$ defines a factor
automorphism on $Q$ which will be denoted by $\ba$. Since $Q$ is
also a $p$-adic Lie group, by \cite{Wa} we get that $U_\pm (\ba )$
are also closed and normal subgroups of $Q$.

Consider $Q/U_-(\ba)$.  Let $\da$ be the factor automorphisms of
$\ba$ defined on $Q/U_-(\ba)$.  It is easy to see that $U_{\pm }(\ap
)C_1(P)\subset U_\pm (\ba)$ and $U_+(\ba )U_-(\ba )\subset U_+(\da
)$: in fact, $U_{\pm }(\ap )C_1(P)= U_\pm (\ba)$ and $U_+(\ba
)U_-(\ba )= U_+(\da )$ follows from \cite{BW}.  Since $P$ is the
closed subgroup generated by $U_-(\ap )$ and $U_+(\ap )$, $U_+(\da )
= Q/U_-(\ba )$. This shows that $Q/U_-(\ba)$ is contracted by $\da$,
that is, $\da ^n (x) \to e$ as $n \to \infty$ uniformly on compact
sets (see \cite{Wa}).  But by assumption $\da$ preserves the Haar
measure on $Q/U_-(\ba)$, hence $Q/U_-(\ba )$ is trivial. This
implies that $P/C_1(P)= U_-(\ba )$. This implies that $P/C_1(P)$ is
contracted by $\ba ^{-1}$. Since $\ba$ preserves the Haar measure on
$P/C_1(P)$, we get that $P/C_1(P)$ is trivial. Thus, $P\subset
C_1(P)$.  Since $P$ is a discrete extension of a solvable group, $P$
is discrete.  Since $P$ is generated b $U_{\pm}$, $P$ is trivial.
\eo

\bo {\bf of Theorem \ref{t1}} Let $\ap$ be the inner-automorphism of
$G$ defined by $g\in G$ and $U_{\pm} = \{ x\in G \mid \lim _{n\to
\pm \infty} \ap ^n (x) = e\}$. Since Ad~$(G)$ is of type $R$, $U_\pm
\subset {\rm Ker~(Ad)}$. Since $\ap ^{\pm 1}|_{U_\pm}$ is a
contraction, there are compact open subgroups $K_{\pm}$ in $U_{\pm}$
such that $\ap ^{-\pm i}(K_\pm) $ is increasing and $U_\pm = \cup
_{i\geq 0} \ap ^{-\pm i}(K_\pm)$.

By assumption $C_k(\ap ^i (K_-), \ap ^{-i}(K_+))$ is discrete.
Since $U_\pm = \cup _{i\geq 0} \ap ^{-\pm i}(K_\pm)$, we get that
$C_k(U_+ , U_-)= \cup _{i\geq 0}C_k(\ap ^i (K_-), \ap ^{-i}(K_+))$
(because an increasing union of closed subgroups is closed in a
$p$-adic Lie group).  Thus, $C_k(U_+, U_-)$ is a countable group,
hence discrete.

By Lemma \ref{l1}, there is a $\ap$-invariant compact subgroup $O$
such that $O$ is centralized by $U_{\pm} $ and $U_+OU_-$ is an open
subgroup of $G$.  Let $N= \langle U_{\pm}, O \rangle$.  Then $N$ is
an $\ap$-invariant open subgroup of $G$.  Considering the group
generated by $N$ and $g$, we get that Haar measures on $N$ and its
quotient are $\ap$-invariant. Since $N$ is generated by $U_{\pm}$
and $O$, $C_k(N/O)$ is discrete. This implies that $N/O$ is a
discrete extension of a solvable group. It follows from Lemma
\ref{topgp} that $U_{\pm}\subset O$.  Since $O$ is compact and
$U_{\pm}$ is a $p$-adic vector space, we get that $U_{\pm }$ is
trivial. This proves that $G$ is of type $R$. \eo

The following provides an example where the condition on the kernel
of the adjoint representation is satisfied.

We say that a finitely generated group $A$ is virtually $\Z ^k$ if
$A$ contains a normal subgroup $B$ of finite index with $B\simeq \Z
^k$.

\bc\label{adj} Let $G$ be a $p$-adic Lie group such that Ad~$(G)$ is
of type $R$ and open subgroups of $G$ and their quotients are
unimodular.  Assume that finitely generated quotients of any open
subgroup of $G$ is virtually $\Z ^k$ for $k\leq 2$ . Then $G$ is of
type $R$. \ec

We first prove the following lemmas.

\bl\label{fincom} If $A$ is a group and $B$ is subgroup of $A$ such
that $B$ is in the center of $A$ and $A/B$ is finite, then $[A, A]$
is finite. \el

\bo Let $a_1, \cdots , a_k$ be such that $A= \cup a_iB$.  Since $B$
is in the center of $A$, for $b_1 ,b_2 \in B$,
$a_ib_1a_jb_2b_1^{-1}a_i^{-1}b_2^{-1} a_j^{-1}=
a_ia_ja_i^{-1}a_j^{-1}$.  Thus, $A$ has finitely many commutators,
hence $[A, A]$ is finite.
\eo

\bl\label{comm} If $A$ is virtually $\Z ^k$ for $k\leq 2$, then
$C_3(A)$ is finite. \el

\bo For any $x\in A$, let $i_x$ denote the innerautomorphism defined
by $x$ on $A$.  If $B$ is a normal subgroup of finite index such
that $B\simeq \Z ^2$, define a homomorphism $f\colon A \to GL(2,
\Z)$ by $f(x) = i_x |_B$.  Since $A/B$ is finite and $B$ is in the
kernel of $f$, $f(A)$ is finite.  Let $A_1 = \Ker (f)$.  Then $A_1$
is a normal subgroup of $A$ such that $A/A_1$ is a 2-step solvable
group - any finite subgroup of $GL(2, \Z)$ is subgroup of some
orthogonal transformations $O(2,\R )$ which is a $2$-step solvable
group.  This shows that $C_2(A)\subset A_1$.

Since $f(x) = i_x|_B$, $B$ is in the center of $A_1$.  Since $A/B$
is finite, $A_1/B$ is finite.  By Lemma \ref{fincom}, $C_1(A_1)$ is
finite.  Since $C_2(A)\subset A_1$, we get that $C_3(A)$ is finite.

If $B \simeq \Z$, define a homomorphism $f'\colon A \to \{ \pm 1 \}$
by $f(x) = i_x |_B$.  Let $A_1 '= \Ker (f')$.  Then $A_1'$ is a
normal subgroup of $A$ such that $A/A_1'$ is an abelian group and
$B$ is in the center of $A_1'$. This shows that $C_1(A)\subset
A_1'$. Since $A_1'/B$ is finite and $B$ is in the center of $A_1'$,
by Lemma \ref{fincom}, we conclude that $C_1(A_1')$ is finite, hence
$C_2(A)$ is finite. \eo

\bo{\bf of Corollary \ref{adj}} Let $F$ be a finite subset of the
kernel of the adjoint representation of $G$ and $H$ be the closed
subgroup generated by $F$.  We now claim that $C_4(H)$ is finite.

Since $F$ is in the kernel of the adjoint representation of $G$,
each $x\in F$ centralizes a compact open subgroup $U_x$ of $G$ (see
Theorem 3, Section 7, Chapter III of \cite{Bo}).  Let $U= \cap U_x$.
Then $U$ is a compact open subgroup of $G$ centralized by all
elements of $F$, hence $U$ is centralized by $H$.

Let $J= UH$.  Then $J/U$ is a finitely generated group.  So, by
assumption $J/U$ is virtually $\Z^k$, for $k\leq 2$.  By Lemma
\ref{comm}, $C_3(J/U)$ is finite.  This implies that $C_3(J)$ is
compact, hence $C_3(H)$ is compact.  Since $U$ is open,
$C_3(H)/C_3(H)\cap U$ is finite.  Since $U$ is centralized by $H$,
by Lemma \ref{fincom} we get that $C_4(H)$ is finite. Now the result
follows from Theorem \ref{t1}\eo

\bigskip\medskip
\advance\baselineskip by 2pt
\begin{tabular}{ll}
C.\ R.\ E.\ Raja \\
Stat-Math Unit \\
Indian Statistical Institute (ISI) \\
8th Mile Mysore Road \\
Bangalore 560 059, India.\\
creraja@isibang.ac.in
\end{tabular}

\end{document}